\documentclass[twoside,reqno,11pt]{fcaa}

\usepackage{graphicx}
\usepackage{epsfig}
\usepackage{amsthm}
\usepackage{amsmath}
\usepackage{latexsym}
\usepackage{amsfonts}
\usepackage{amssymb}

\newtheoremstyle{theorem}
  {15pt}          
  {15pt}  
  {\sl}  
  {\parindent}
  {\sc}  
  {. }    
  { }    
  {}     
\theoremstyle{theorem}

\newtheorem{theorem}{Theorem}[section]
\newtheorem{corollary}{Corollary}[section]

\newtheoremstyle{defi}
  {15pt}          
  {15pt}  
  {\rm}  
  {\parindent}     
  {\sc}  
  {. }    
  { }    
  {}     
\theoremstyle{defi}

\newtheorem{remark}{Remark}[section]
\newtheorem{example}{Example}[section]

 \def\proofname{\indent\rm P\ r\ o\ o\ f.}
 
 \def\qed{\hfill$\Box$} 

\usepackage{hyperref} 
\usepackage{amsfonts}

 \def\theequation{\arabic{section}.\arabic{equation}}

 \def\themycopyrightfootnote{\vspace*{3pt}
   October 12, 2015
  }

 \title[MONOTONICITY AND SIGN OF DERIVATIVE \dots]{MONOTONICITY OF FUNCTIONS\\ [4pt] AND SIGN CHANGES\\ [4pt] OF THEIR CAPUTO DERIVATIVES}
 \author[K. Diethelm]{Kai Diethelm$^{1,2}$}

\begin{document}

 \vbox to 2.5cm { \vfill }


 \bigskip \medskip

\begin{abstract}

  It is well known that a continuously differentiable function is monotone in an
  interval $[a,b]$ if and only if its first derivative does not change its sign
  there. We prove that this is equivalent to requiring that the Caputo
  derivatives of all orders $\alpha \in (0,1)$ with starting point $a$
  of this function do not have a change of sign there. In contrast to
  what is occasionally conjectured, it not sufficient if the Caputo
  derivatives have a constant sign for a few values of $\alpha \in
  (0,1)$ only.

  \medskip


  \smallskip


\end{abstract}

\maketitle

\vspace*{-16pt}


\section{Introduction}\label{sec:intro}

\setcounter{section}{1}
\setcounter{equation}{0}\setcounter{theorem}{0}

Classical integer order calculus provides numerous results that
link the shape properties of a function's graph (such as,
e.\ g., monotonicity or the location of extrema) with the signs of its
first or second 
derivative. The following theorems contain some specific
examples; since they are so well known, we omit the proofs.

\begin{theorem}\label{thm:c-mon}
  A function $f \in C^1[a,b]$ is monotone on $[a,b]$ if and only if
  its first derivative $f'$ does not change its sign in this interval.
\end{theorem}

\begin{theorem}\label{thm:c-extr}
  If a function $f \in C^1[a,b]$ attains a global maximum or a global
  minimum at some point $t^* \in (a,b)$ then $f'(t^*) = 0$.
  
  If a function $f \in C^2[a,b]$ attains a global maximum 
  at some point $t^* \in (a,b)$ then $f'(t^*) = 0$ and $f''(t^*)
  \le 0$; if $t^* \in (a,b)$ is the location of a global minimum of
  $f$ then $f'(t^*) = 0$ and $f''(t^*) \ge 0$.
\end{theorem}

It is quite natural to ask for generalizations of such statements that
involve fractional, rather than integer order,
derivatives. In this paper we concentrate on fractional derivatives
of Caputo type \cite[Chapter 3]{Diethelm}, i.\ e.\  
\begin{equation}
  \label{eq:caputo}
  D^\alpha_{*a} f(t) := J^{\lceil \alpha \rceil - \alpha}_a D^{\lceil \alpha \rceil} f(t)
\end{equation}
for $\alpha > 0$ and $\alpha \notin \mathbb N$
where $\lceil \cdot \rceil$ denotes the ceiling function, $D^m$ is the
standard $m$th order differential operator for $m \in \mathbb N$ and
$J^\beta_a$ is the Riemann-Liouville integral operator of order
$\beta>0$ with starting point $a$, viz.\ 
\begin{equation}
  \label{eq:rli}
  J^\beta_a g(t) := \frac1{\Gamma(\beta)} \int_a^t (t-s)^{\beta-1} g(s) {\, \mathrm d} s.
\end{equation}
In this case, the following fractional extension of the first part of
Theorem \ref{thm:c-extr} has been shown by Luchko \cite[Theorem 1]{Luchko}:

\begin{theorem}\label{thm:f-extr-1}
  If a function $f \in C^1[a,b]$ attains a global maximum 
  at some point $t^* \in (a,b)$ then $D^\alpha_{*a} f(t^*) \ge 0$ for
  all $\alpha \in (0,1)$. If $t^* \in (a,b)$ is the location of the
  global minimum of $f$ then $D^\alpha_{*a} f(t^*) \le 0$ for
  all $\alpha \in (0,1)$.  
\end{theorem}

In fact, Luchko's result is more general than our formulation
here since he admits functions $f$ from a slightly larger class of
functions; however, for the sake of simplicity, we have chosen this
marginally restricted version.

Some authors have believed that this result can be extended to
fractional derivatives of order $\alpha \in (1,2)$; specifically they
claimed that, under the additional condition $f \in C^2[a,b]$, the
inequality 
$$ D^\alpha_{*a} f(t^*) \ge 0 $$
holds for these $\alpha$ if $t^*$ is a global minimum point of
$f$. While Luchko's result (Theorem \ref{thm:f-extr-1}) extends the
first part of the classical Theorem \ref{thm:c-extr}, such a result
would provide an extension of the  second part of Theorem
\ref{thm:c-extr}. However, recently Al-Refai \cite{Al-Refai}
(who also provides some references to claims of such results in the
literature) explicitly constructed a counterexample and hence
demonstrated that the second part of Theorem \ref{thm:c-extr} cannot
be extended in this manner. 

Our goal here is to perform a similar investigation for possible fractional versions
of Theorem \ref{thm:c-mon}. In a certain sense, it thus continues the
author's earlier work \cite{Diethelm-mvt} that transferred other
classical calculus results, including the mean value theorems of differential
and integral calculus, to the fractional setting.

\section{Caputo derivatives of monotone functions}\label{sec:main}

\setcounter{section}{2}
\setcounter{equation}{0}\setcounter{theorem}{0}

We begin our investigation of a generalization of Theorem
\ref{thm:c-mon} to the fractional case with a very simple but
nevertheless useful observation.

\begin{theorem} \label{thm1}
  Let the function $f \in C^1[a,b]$ be monotone on $[a,b]$. Then, for
  all $\alpha \in (0,1)$ and all $t \in [a,b]$, we have that
  $$
     D^\alpha_{*a} f(t) 
      \begin{cases}
           \le 0 & \mbox{if } f \mbox{ is decreasing,}  \\
           \ge 0 & \mbox{if } f \mbox{ is increasing.}
      \end{cases}
  $$
\end{theorem}

In this formulation, we use the terms ``increasing'' and
``decreasing'' in the non-strict sense.

\proof
  Assume that $f$ is increasing. Then, $f'$ is nonnegative in $[a,b]$,
  and hence, by the definition of the Caputo operator given in
  eq.\ (\ref{eq:caputo}) and eq.\ (\ref{eq:rli}),
  $$ 
    D^\alpha_{*a} f(t) 
     = \frac1{\Gamma(1-\alpha)} \int_a^t (t-s)^{-\alpha} f'(s) \, {\mathrm d} s
     \ge 0.
  $$
  In the case of a decreasing function $f$, we argue in an analog way.
\endproof

We obtain an immediate consequence that can be interpreted as the fractional
generalization of the ``forward direction'' of Theorem
\ref{thm:c-mon}.

\begin{corollary}
  Let the function $f \in C^1[a,b]$ be monotone on $[a,b]$. Then, for
  all $\alpha \in (0,1)$, the Caputo derivative $D^\alpha_{*a} f$ does
  not have a change of sign in $[a,b]$.  
\end{corollary}

But we can also obtain a fractional generalization of the ``backward
direction'' of Theorem \ref{thm:c-mon}: 

\begin{theorem} \label{thm2}
  Assume that $f \in C^1[a,b]$ is such that $D^\alpha_{*a} f(t) \ge 0$
  for all $t \in [a,b]$ and all $\alpha \in (\alpha_0, 1)$ with some
  $\alpha_0 \in (0,1)$. Then, $f$ is monotone increasing.
  Similarly, if $D^\alpha_{*a} f(t) \le 0$ for all $t$ and $\alpha$
  mentioned above, then $f$ is monotone decreasing.
\end{theorem}

\proof
  Our observation in the first case is based on the definition of the
  Caputo derivative and certain well known properties of the
  Riemann-Liouville integral. Specifically, consider first the case that $t \in
  (a, b]$. Then, in view of eq.\ (\ref{eq:caputo}) and our assumption, 
  $0 \le D^\alpha_{*a} f(t) = J^{1-\alpha}_a f'(t)$. Since $f'$ is
  assumed to be continuous, we may invoke \cite[Theorem
    2.10]{Diethelm} to determine the limit $\alpha \to 1-$
  of the right-hand side of this inequality; this yields
  \begin{eqnarray*}
     0 &\le & \lim_{\alpha \to 1-} D^\alpha_{*a} f(t)
          = \lim_{\alpha \to 1-} J^{1-\alpha}_a f'(t)
          = f'(t).
  \end{eqnarray*}
  We thus conclude that $f'$ is nonnegative on the half-open interval
  $(a,b]$. The continuity of $f'$ on the closed interval $[a,b]$ then
  implies that $f'(a) \ge 0$ as well; hence $f'$ is nonnegative on the
  entire interval $[a,b]$, and it follows that $f$ is increasing.

  The second claim follows immediately if we apply the result that we
  have just shown to $-f$ instead of $f$.
\endproof

Theorem \ref{thm2} has an interesting consequence that we explicitly
note.

\begin{theorem}
  Assume that $f \in C^1[a,b]$ is such that $D^\alpha_{*a} f(t) \ge 0$
  for all $t \in [a,b]$ and all $\alpha \in (\alpha_0, 1)$ with some
  $\alpha_0 \in (0,1)$. Then, $D^\alpha_{*a} f(t) \ge 0$
  for all $t \in [a,b]$ and all $\alpha \in (0, 1)$. 

  The analog statement holds when the sign of the derivatives is
  reversed. 
\end{theorem}

Thus, if we can prove that the Caputo derivatives of $f$ are of a
constant sign for all $\alpha$ sufficiently close to $1$, then it
already follows that this property even holds for all $\alpha$ from the
complete interval $(0,1)$ and, by Theorem \ref{thm2}, also for the
classical case $\alpha = 1$.

\proof
  If the Caputo derivatives $D^\alpha_{*a} f(t)$ are nonnegative for
  the indicated values of $t$ and $\alpha$ then Theorem \ref{thm2}
  tells us that $f$ is increasing; the claim then follows from Theorem
  \ref{thm1}. The statement for the case of opposite signs of
  $D^\alpha_{*a} f(t)$ follows in an analog manner.
\endproof

It thus turns out that, contrary to an occasionally stated conjecture,
the lack of a change of sign of $D^\alpha_{*a} f$ for
some, but not all, $\alpha \in (0,1)$ is not sufficient to imply the
monotonicity of $f$. One may ask how a non-monotone functions which
has a fractional derivative without a change of sign may look like,
and in particular whether such a function must necessarily be of a
very exotic or even pathological nature. The following example
demonstrates that this is not the case; indeed a simple polynomial of
a quite low degree can already have this property. An interesting
aspect in this connection is the fact that our 
counterexample is the same one that Al-Refai used in his
investigations from \cite[Section 2]{Al-Refai} regarding the
fractional version of the second part of 
Theorem \ref{thm:c-extr} that we mentioned above.

\begin{example}\label{Ex1}
  Define the function $g$ by 
  $$
     g(t) := t \left( t - \frac12\right) (t-1) = t^3 - \frac32 t^2 + \frac12 t,
  $$
  so that
  \begin{equation}
    \label{eq:g1}
    g'(t) = 3 t^2 - 3 t + \frac 12.
  \end{equation}
  and
  \begin{eqnarray}
     D^\alpha_{*0} g(t) &=&  \nonumber
          \frac6{\Gamma(4-\alpha)} t^{3-\alpha} - \frac3{\Gamma(3-\alpha)} t^{2-\alpha} + \frac1{2\Gamma(2-\alpha)} t^{1-\alpha} \\
      &=& \frac{6t^{1-\alpha}}{\Gamma(4-\alpha)} \left(t^2 - \frac12 (3-\alpha) t + \frac1{12}(3-\alpha)(2-\alpha)\right)
          \label{eq:ga}
  \end{eqnarray}
  for $0 < \alpha < 1$.
  It follows from eq.\ (\ref{eq:ga}) that the only zeros of
  $D^\alpha_{*0} g$ are located at $t_1(\alpha) = 0$, 
  $t_2(\alpha) = \left(3 - \alpha - \sqrt{(3-\alpha) (1+\alpha) / 3}\right) / 4$
  and
  $t_3(\alpha) = \left(3 - \alpha + \sqrt{(3-\alpha) (1+\alpha) / 3}\right) / 4$. 
  Clearly, $t_1(\alpha) < t_2(\alpha) \le t_3(\alpha)$, 
  and hence we conclude that $D^\alpha_{*0} g$
  does not have a change of sign on $[0, t_2(\alpha)]$.

  But now, 
  $$
    \frac{\mathrm d}{\mathrm d \alpha} t_2(\alpha)
     = \frac14 \left(-1 - \frac{1-\alpha}{\sqrt{3 (3-\alpha) (1+\alpha)}} \right)
     < 0
  $$
  for $0 < \alpha < 1$. Therefore, if we denote by $I_\alpha :=
  [0, t_2(\alpha)]$ the largest interval that has the left end point at $0$
  and that is such that $D^\alpha_{*0} g$ does not have a change of
  sign on $I_\alpha$, then these intervals monotonically shrink as
  $\alpha$ increases. Specifically,
  $$
    \frac12 = t_2(0) > t_2(\alpha_1) > t_2(\alpha_2) > t_2(1) = \frac{3-\sqrt3}6  \approx 0.211325,
  $$
  for $0 < \alpha_1 < \alpha_2 < 1$, and hence 
  $$ 
     \left[0, \frac12\right] = I_0 \supset I_{\alpha_1} \supset I_{\alpha_2} 
                             \supset I_1 = \left[ 0, \frac{3-\sqrt3}6 \right].
  $$

  Therefore, if we pick some $\alpha^* \in (0,1)$, then we can see
  that $D^{\alpha^*}_{*0} g$ is nonnegative throughout the interval
  $I^* := I_{\alpha^*} = [0, t_2(\alpha^*)]$. Moreover, it follows from our
  monotonicity observation on $t_2$ that all derivatives
  $D^{\alpha}_{*0} g$ with $\alpha \in (0,\alpha^*)$ are nonnegative
  on the interval $I^*$ as well. Nevertheless, an elementary
  computation reveals that the only zeros of $g'$ are located at 
  $\tau_1 = (3 - \sqrt3)/6 = t_2(1)$ and $\tau_2 = (3 + \sqrt3)/6 =
  1 - \tau_1 \approx 0.788675$. From this observation and the fact that
  $g$ is a polynomial of degree 3 with a positive leading coefficient,
  we conclude that $g$ has a local maximum at the point $t_2(1)$ which
  is located in the interior of our interval $I^*$; hence, in spite of
  the fact that $D^{\alpha}_{*0} g(t) > 0$ for almost all $t \in I^*$
  (with the sole exception of the point $t=0$ for which
  $D^{\alpha}_{*0} g(t) =0$ holds) and all $\alpha \in (0, \alpha^*]$, 
  the function $g$ is not monotone on $I^*$.
\end{example}

\begin{remark}
  We point out that all our investigations are based on the fact that
  the left end point of the interval on which the monotonicity is
  discussed coincides with the starting point of the Caputo
  differential operator. This is a very natural requirement because of
  the non-locality of the fractional differential
  operator. Specifically, if we had $a < a' < b$ and wanted to
  find out whether a function $f$ is monotonic on $[a',b]$, then a
  criterion based on the behaviour of $D^\alpha_{*a} f$ would not be of
  any help since the values of $D^\alpha_{*a} f$ may de influenced to
  an arbitrarily large extent by the behaviour of $f$ on the interval
  $[a, a']$, i.\ e.\ outside of the interval in which we are
  interested. Therefore, in a situation like this only integer-order
  derivatives can be of any use because they are local operators and
  hence do not depend on the behaviour of $f$ outside of the interval
  of interest.
\end{remark}




 \bigskip \smallskip

 \it

 \noindent
$^1$ AG Numerik\\
Institut Computational Mathematics\\
Technische Universit\"at Braunschweig\\
Pockelsstr.\ 14\\
38106 Braunschweig, GERMANY\\
e-mail: k.diethelm@tu-bs.de
 \\[12pt]
$^2$ GNS Gesellschaft f\"ur numerische Simulation mbH\\
Am Gau\ss berg 2\\
38114 Braunschweig, GERMANY\\
e-mail: diethelm@gns-mbh.com

\end{document}